\documentclass[11pt,twoside]{article}

\usepackage{amsfonts}
\usepackage[mathscr]{euscript}
\usepackage{xcolor}

\newcommand{\lbl}[1]{\label{#1}}

\newtheorem{theo}{Theorem}[section]
\newtheorem{prop}{Proposition}[section]
\newtheorem{lem}{Lemma}[section]
\newtheorem{remark}{Remark}[section]
\newtheorem{cor}{Corollary}[section]

\newcommand{\be}{\begin{equation}}
\newcommand{\ee}{\end{equation}}
\newcommand\bes{\begin{eqnarray}} \newcommand\ees{\end{eqnarray}}
\newcommand{\bess}{\begin{eqnarray*}}
\newcommand{\eess}{\end{eqnarray*}}
\newcommand\bedd{\bes\left\{\begin{array}{ll}\smallskip}
\newcommand\eedd{\end{array}\right.\ees}

\begin{document}\thispagestyle{empty}

\begin{center}{\Large\bf A diffusive logistic equation with a free boundary and }\\[2mm]
  {\Large\bf sign-changing coefficient in time-periodic environment}\footnote{This work was supported by NSFC Grant 11371113}\\[5mm]
 {\Large  Mingxin Wang\footnote{E-mail: mxwang@hit.edu.cn; Tel: 86-15145101503; Fax: 86-451-86402528}}\\[2mm]
{Natural Science Research Center, Harbin Institute of Technology, Harbin 150080, PR China}
\end{center}

\begin{quote}
\noindent{\bf Abstract.} This paper concerns a diffusive logistic equation with a free boundary and sign-changing intrinsic growth rate in heterogeneous time-periodic environment, in which the variable intrinsic growth rate may be ``very negative" in a ``suitable large region" (see conditions {\bf(H1)}, {\bf(H2)}, (\ref{4.3})). Such a model can be used to describe the spreading of a new or invasive species, with the free boundary representing the expanding front. In the case of higher space dimensions with radial symmetry and the intrinsic growth rate has a positive lower bound, this problem has been studied by Du, Guo \& Peng \cite{DGP}. They established a spreading-vanishing dichotomy, the sharp criteria for spreading and vanishing and estimate of the asymptotic spreading speed. In the present paper, we show that the above results are retained for our problem.


 \end{quote}

\def\theequation{\arabic{section}.\arabic{equation}}

 \section{Introduction}
 \setcounter{equation}{0} {\setlength\arraycolsep{2pt}

In the real world, the invasion of harmful species and/or introduction of new and beneficial species are natural phenomena. It is of primary importance to predict and analyze the growing and spreading mechanism of biological invasions. A lot of mathematicians have made efforts to develop various invasion models and investigated them from a viewpoint of mathematical ecology. Most theoretical approaches are based on or start with single-species models. In consideration of the heterogeneous environment, the following problem
\bess
 \left\{\begin{array}{lll}
 u_t-d\Delta u=a(t,x)u-b(t,x)u^2, \ \ &t>0,\ \ x\in\Omega,\\[1mm]
 B[u]=0,  &t\ge 0, \ \ x\in\partial\Omega,\\[1mm]
 u(0,x)=u_0(x),  & x\in\Omega
 \end{array}\right.
 \eess
is a typical one to describe the spread, persistence and extinction of the new or invasive species and has received an astonishing amount of attention, please refer to \cite{CC1}-\cite{CLL} and \cite{DP, HC, Hess, LL, Ni} for example. In this model, $u(t,x)$ represents the population density; constant $d>0$ denotes the diffusion (dispersal) rate; $a(t,x)$ and $b(t,x)$ represent the intrinsic growth rate and self-limitation coefficient of the species, respectively; $\Omega$ is a bounded domain of $\mathbb R^N$; the boundary operator $B[u]=\alpha u+\beta \frac{\partial u}{\partial \nu}$, $\alpha$ and $\beta$ are non-negative functions and satisfy $\alpha+\beta>0$, $\nu$ is the outward unit normal vector of the boundary $\partial\Omega$. The corresponding systems with heterogeneous environment have also been studied extensively, please refer to \cite{CC2, CLL, HMP, LL, Ni} and the references cited therein.

In most spreading processes in the natural world, a spreading front can be observed. When a new or invasive species initially occupies a region $\Omega_0$ with density $u_0(x)$, as time $t$ increases, it is natural to expect that $\Omega_0$ will evolve into an expanding region $\Omega(t)$ with an expanding front $\partial\Omega(t)$, inside which the initial function $u_0(x)$ will evolve into a positive function $u(t,x)$ governed by a suitable diffusive equation, with $u(t,x)$ vanishing on the moving boundary $\partial\Omega(t)$.

In the natural world, for most animals and plants, their birth and death rates will change with the seasons, so the intrinsic growth rate $a(t,x)$ and then the self-limitation coefficient $b(t,x)$ should be time-periodic functions. Especially, in the winter of severe cold and cold zones, animals cannot capture enough food to feed
upon and do not breed, seeds cannot germinate and buds cannot grow above ground, so their birth rates are zero. In the meantime, their death rates will be greater. Therefore, in some periods and some areas, the intrinsic growth rate $a(t,x)$ may be negative. In order to simplify the mathematics, in this paper we only consider the one dimensional case, i.e., $N=1$, and assume that the left boundary is fixed: $x=0$. For a more realistic description of the growth mechanism and spreading of a new or invasive species, throughout this paper we assume
\vspace{-1mm}\begin{quote}
 {\bf(H)}\, Functions $a,\,b\in \big(C^{\nu/2,\,\nu}\cap
 L^\infty\big)([0,\infty)\times[0,\infty))$ for some $\nu\in(0,1)$, and are $T$-periodic in time $t$ for some $T>0$. The function $a(t,x)$ is positive somewhere in $[0,T]\times[0,\infty)$, the function $b(t,x)$ satisfies $b_1\le b(t,x)\le b_2$ in $[0,\infty)\times[0,\infty)$ for some  positive constants $b_1$, $b_2$.
\vspace{-1mm}\end{quote}

Motivated by the above natural phenomena, in this paper we shall study the following free boundary problem
 \bes
 \left\{\begin{array}{lll}
 u_t-du_{xx}=a(t,x)u-b(t,x)u^2, \ \ &t>0,\ \ 0<x<h(t),\\[1mm]
 B[u](t,0)=0,\ \ u(t,h(t))=0,\ \ &t\ge 0,\\[1mm]
 h'(t)=-\mu u_x(t,h(t)),&t\ge 0,\\[1mm]
 h(0)=h_0,\ \ u(0,x)=u_0(x),&0\le x\le h_0,
 \end{array}\right.\label{1.1}
 \ees
where, $h_0$ denotes the size of initial habitat, $\mu$ is the ratio of expanding speed of the free boundary and population gradient at expanding front, it can also be considered as the ``moving parameter" of the free boundary, $B[u]=\alpha u-\beta u_x$, $\alpha$ and $\beta$ are non-negative constants and satisfy $\alpha+\beta=1$, the initial function $u_0(x)$ satisfies
 \vspace{-1mm}\begin{quote}
  $u_0\in C^2([0,h_0])$ , $u_0>0$ in $(0,h_0)$, $B[u_0](0)=u_0(h_0)=0$.
  \vspace{-1mm}\end{quote}
The free boundary condition $h'(t)=-\mu u_x(t,h(t))$ in (\ref{1.1}) is a one-phase Stefan condition and can refer to \cite{BDK} for the more detailed background.

When both functions $a$ and $b$ have positive lower and upper bounds, i.e., there exist positive constants $\kappa_1$, $\kappa_2$ such that
 \[\kappa_1\le a(t,x),\,b(t,x)\le \kappa_2,\]
recently, Du, Guo \& Peng in \cite{DGP} have studied the problem (\ref{1.1}) for the higher dimension and radially symmetric case. By developing the earlier techniques and introducing new ideas and methods, they have obtained various interesting results. When $b=1$ and $a=a(x)$ changes sign, the problem (\ref{1.1}) has been studied by Zhou \& Xiao \cite{ZX} and Wang \cite{Wang14}. Complete conclusions about spreading-vanishing dichotomy, sharp criteria for spreading and vanishing and asymptotic spreading speed of the free boundary were obtained by the author.

The main aim of this paper is to study the dynamics of (\ref{1.1}) and show that most results of \cite{DGP, Wang14, ZX} continue to hold in the more realistic situation that the $T$-periodic intrinsic growth rate $a(t,x)$ may be ``very negative" in the sense that both $|\{a(t,x)>0\}|\ll|\{a(t,x)<0\}|$ and $\int_0^\infty a(t,x){\rm d}x=-\infty$ are allowed for any $t>0$. To our best knowledge, the present paper seems to be the first attempt to consider the sign-changing intrinsic growth rate with the time-periodic and spatial heterogeneous environment in the moving domain problem.

In the special case that the functions $a$ and $b$ are independent of time $t$ and have positive lower and upper bounds, such kind of problems has been studied systematically. For example, when $a$ and $b$ are positive constants, the problem (\ref{1.1}) was investigated earlier by \cite{DLin} for $\alpha=0$ and by \cite{KY} for $\beta=0$; when $a=a(x)$ and $b=b(x)$, the problem (\ref{1.1}) was discussed by Du, Guo \& Liang \cite{DG, DLiang} for the higher dimension and radially symmetric case;  the non-radial case in higher dimensions was treated by Du \& Guo \cite{DG1}. Instead of $u(a-bu)$ by a general function $f(u)$, Du, Matsuzawa \& Zhou \cite{DMZ}, Kaneko \cite{Ka} and Du \& Lou \cite{DLou} investigated the corresponding free boundary problems.

Peng \& Zhao \cite{PZ} studied a free boundary problem of the diffusive logistic model with seasonal succession. They considered that the species does not migrate and stays in a hibernating status in the bad season. The evolution of the species obeys the Malthusian equation $u_t=-\delta u$ in the bad season, and obeys the diffusive logistic equation with positive constant coefficients in the good season. The diffusive competition system with positive constant coefficients and a free boundary has been studied by Guo \& Wu \cite{GW}, Du \& Lin \cite{DL2} and Wang \& Zhao \cite{WZjdde}. The diffusive prey-predator model with positive constant coefficients and free boundaries has been studied by Wang \& Zhao \cite{Wjde, WZ, ZW}.

This paper is organized as follows.  In Section 2, we prove the global existence, uniqueness, regularity and estimate of $(u,h)$. Especially, the uniform estimates of $\|u(t,\cdot)\|_{C^1[0,\,h(t)]}$ for $t\ge 1$ and $\|h'\|_{C^{(1+\nu)/2}([n+1, n+3])}$ for $n\ge 0$ are obtained directly regardless of the size of $h_\infty:=\lim_{t\to\infty}h(t)$. These uniform estimates allow us to assert $h'(t)\to 0$ when $h_\infty<\infty$ and play a key role for determining the vanishing phenomenon. In Section 3, we exhibit some fundamental results, including the comparison principle for the moving domain, and some properties of the principal eigenvalue of a $T$-periodic eigenvalue problem. In Section 4, we shall derive the spreading-vanishing dichotomy:

Either
\vspace{-2mm}\begin{itemize}
\item[(i)]
Spreading: $\lim_{t\to\infty}h(t)=\infty$ and $\lim_{n\to\infty}u(t+nT,x)=U(t,x)$ uniformly on $[0,T]\times[0,L]$ for any $L>0$, where $U(t,x)$ is the unique $T$-periodic positive solution to {\rm(\ref{4.2})};\end{itemize}\vspace{-2mm}
or
\begin{itemize}
\vspace{-2mm}\item[(ii)] Vanishing: $\lim_{t\to\infty}h(t)=h_\infty<\infty$ and
 $\lim_{t\to\infty}\,\max_{0\leq x\leq h(t)}u(t,x)=0$.
 \end{itemize}\vspace{-1mm}
In Section 5, we first establish the sharp criteria for spreading and vanishing, then estimate the asymptotic spreading speed of the free boundary. The last section is a brief discussion.

We should remark that for the higher dimensional and radially symmetric case of (\ref{1.1}), the methods in this paper are still valid and the corresponding results can be retained.

\section{Global existence, uniqueness, regularity and estimate of solution}
\setcounter{equation}{0}{\setlength\arraycolsep{2pt}

In this section, we give the existence, uniqueness and estimate of the solution $(u,h)$ to (\ref{1.1}).

\begin{theo}\lbl{th2.1}\, The problem {\rm(\ref{1.1})} has a unique global solution $(u,h)$ and
 \bes
  u\in C^{1+\frac{\nu} 2,2+\nu}(D_\infty),\ \ h\in C^{1+\frac{1+\nu}2}((0,\infty)),
  \lbl{2.1}\ees
where
 \[D_\infty=\big\{(t,x):\,t>0,\, x\in\big[0,h(t)\big]\big\}.\]
Moreover, there exists a positive constant $M=M\left(\|a,b,u_0\|_\infty\right)$ such that
  \bes
  0<u(t,x)\leq M, \ \ 0<h'(t)\leq \mu M, \ \ \forall \  t> 0,\ 0<x<h(t).\lbl{2.2}\ees
Further more, if $h_\infty:=\lim_{t\to\infty}h(t)<\infty$, then there exists a positive constant $C=C\left(\mu,h_\infty,M\right)$ such that
   \bes&\|h'\|_{C^{\nu/2}([n+1,n+3])}\leq C, \ \ \forall \  n\geq 0,&\lbl{2.3}\\[1mm]
  &\|u(t,\cdot)\|_{C^1([0,\,h(t)])}\leq C, \ \ \forall \  t\ge 1.&
 \lbl{2.4}\ees
\end{theo}

{\bf Proof.}\, Noting that the functions $a(t,x)$ and $b(t,x)$ are bounded, applying the methods used in \cite{ChenA, DLin} with some modifications, we can prove that (\ref{1.1}) has a unique global solution $(u,h)$, $u\in C^{\frac{1+\nu} 2,1+\nu}(D_\infty)$, $h\in C^{1+\frac{\nu}2}(0,\infty)$ and $u$ satisfies the first estimate of (\ref{2.2}). The second estimate of (\ref{2.2}) can be proved by the similar way to that of \cite[Theorem 2.1]{Wang14}. The details are omitted here.

\vskip 2pt Now, let us prove (\ref{2.1}). Let $y=x/h(t)$ and $w(t,y)=u(t,x)$. A simple calculation gives
 \bes
\left\{\begin{array}{ll}
w_t-d\zeta(t)w_{yy}-\xi(t,y)w_y=a\big(t,h(t)y\big)w-b\big(t,h(t)y\big)w^2, \ &t>0,\ 0<y<1,\\[2mm]
\left(\alpha w-\frac{\beta}{h(t)}w_y\right)(t,0)=0, \ \ w(t,1)=0,\ \ &t\ge 0,\\[2mm]
w(0,y)=u_0(h_0y),&0\leq y\leq 1,
\end{array}\right.
 \lbl{2.5}\ees
where $\zeta(t)=h^{-2}(t)$, $\xi(t,y)=yh'(t)/h(t)$. This is an initial-boundary value problem with fixed boundary. Remember $a,\,b\in C^{\frac{\nu}2,\nu}([0,\infty)\times[0,\infty))$, $u\in C^{\frac{1+\nu} 2,1+\nu}(D_\infty)$ and $h\in C^{1+\frac{\nu}2}(0,\infty)$. For any given $\tau>0$ and $0<\varepsilon\ll 1$, applying Theorem 10.1 of \cite[Chap.4, p.351]{LSU} to the problem (\ref{2.5}) in $[\varepsilon,\tau]\times[\varepsilon,1]$ and $[\varepsilon,\tau]\times[0,1-\varepsilon]$, respectively, we obtain that
  \[w\in C^{1+\frac{\nu}2,\,2+\nu}\big([\varepsilon,\tau]\times[\varepsilon,1]\big)
  \cap C^{1+\frac{\nu}2,\,2+\nu}\big([\varepsilon,\tau]\times[0,1-\varepsilon]\big).\]
This implies
 \bess u\in C^{1+\frac{\nu}2,\,2+\nu}\big([\varepsilon,\tau]\times[\varepsilon h(t),h(t)]\big)
  \cap C^{1+\frac{\nu}2,\,2+\nu}\big([\varepsilon,\tau]\times[0,(1-\varepsilon)h(t)]\big).
  \eess
Due to the arbitrariness of $\varepsilon$, one achieves
 \[u\in C^{1+\frac{\nu}2,\,2+\nu}(D_\tau),  \ \ \Longrightarrow \ u_x\in C^{\frac{1+\nu}2,\,1+\nu}(D_\tau),\]
where $D_\tau=\{(t,x): 0<t\le\tau,\, 0\le x\le h(t)\}$. Hence, by the condition $h'(t)=-\mu u_x(t,h(t))$, it is immediately to get $h'\in C^{\frac{1+\nu}2}((0,\tau])$.

Finally, we prove (\ref{2.3}) and (\ref{2.4}). For the integer $n\ge 0$, let $w^n(t,y)=w(t+n,y)$. Then $w^n$ satisfies
 \bess
\left\{\begin{array}{ll}
w^n_t-d\zeta^n w^n_{yy}-\xi^n w^n_y=f^n(t,y), \ &0<t\le 3,\ 0<y<1,\\[2mm]
\left(\alpha w^n-\frac{\beta}{h^n}w^n_y\right)(t,0)=0, \ \ w^n(t,1)=0,\ \ &0\le t\le 3,\\[2mm]
w^n(0,y)=u(n,h(n)y),&0\leq y\leq 1,
\end{array}\right.
 \eess
where $\zeta^n=\zeta(t+n)$, $\xi^n=\xi(t+n,y)$, $h^n=h(t+n)$, and
 \[f^n(t,y)=a(t+n,h^n(t)y)w^n(t,y)-b(t+n,h^n(t)y)(w^n(t,y))^2.\]
Noticing (\ref{2.2}), we have that $w^n$, $\zeta^n$, $\xi^n$ and $f^n$ are bounded uniformly on $n$, and
 \[\omega^n(r)=\max_{0\le s,t\le 3,\, |s-t|\le r}|\zeta^n(s)-\zeta^n(t)|\le \frac{2}{h_0^3}\mu Mr\to 0 \ \ \mbox{as} \ \ r\to 0\]
uniformly on $n$. Moreover, $\zeta^n(t)\ge h_\infty^{-2}$ for all $n\ge 0$ and $0\le t\le 3$ as  $h(t)\le h_\infty<\infty$.

Remember the boundary conditions of $w^n(t,y)$ at $y=0,1$ and the fact $0<(h^n(t))'\le\mu M$. Choose $p\gg 1$, we can apply the interior $L^p$ estimate to derive that there exists a positive constant $C$ independent of $n$ such that, for all $n\ge 0$,
 \vspace{-1mm}\begin{quote}
(i)\, when $\beta=0$, we have $\|w^n\|_{W^{1,2}_p([1,3]\times[0,1])}\le C$ (see \cite[Theorem 7.15]{Lie});

(ii)\, when $\beta>0$, we have $\|w^n\|_{W^{1,2}_p([1,3]\times[\frac 12,1])}\le C$ (see \cite[Theorem 7.15]{Lie}), $\|w^n\|_{W^{1,2}_p([1,3]\times[0,\frac 12])}\\ \le C$ (see \cite[Theorem 7.20]{Lie}).
\vspace{-1mm}\end{quote}
In a word, $\|w^n\|_{W^{1,2}_p([1,3]\times[0,1])}\le C$ for all $n\ge 0$. In view of the embedding theorem, it follows that $\|w^n\|_{C^{\frac{1+\nu}2, 1+\nu}([1,3]\times[0,1])}\le C$ for all $n\ge 0$. This implies $\|w\|_{C^{\frac{1+\nu}2, 1+\nu}(E_n)}\le C$ for all $n\ge 0$, where $E_n=[n+1,n+3]\times[0,1]$. This fact combined with
 \[h'(t)=-\mu u_x(t,h(t)), \ u_x(t,h(t))=h^{-1}(t)w_y(t,1), \ 0<h'(t)\le\mu M,\]
allows us to derive (\ref{2.3}). Since these rectangles $E_n$ overlap and $C$ is independent of $n$, it follows that
$\|w\|_{C^{0, 1}([1,\infty)\times[0,1])}\le C$. Using $u_x=h^{-1}(t)w_y$ again, we get (\ref{2.4}). The proof is complete.\ \ \ \fbox{}

\section{Preliminaries}
\setcounter{equation}{0}{\setlength\arraycolsep{2pt}

In this section, we first state a comparison principle and then show some properties of the principal eigenvalue of a $T$-periodic eigenvalue problem. Finally, we discuss the existence and uniqueness of positive solution to a $T$-periodic boundary value problem in the bounded interval.

The following lemma is the analogue of Lemma 3.5 in \cite{DLin} and the proof
will be omitted.

\begin{lem}$($Comparison principle$)$\label{l3.1}\, Let $\overline h\in C^1([0,\infty))$ and $\overline h>0$ in $[0,\infty)$, $\overline u\in C^{0,1}(\overline{Q})\cap C^{1,2}(Q)$, with $Q=\{(t,x): t>0,\, 0<x<\overline h(t)\}$. Assume that $(\overline u, \overline h)$ satisfies
 \bess\left\{\begin{array}{ll}
  \overline u_t-d\overline u_{xx}\geq a(t,x)\overline u-b(t,x)\overline u^2,\ \ &t>0,\ \ 0<x<\overline h(t),\\[1mm]
 B[\overline u](t,0)\geq 0,\ \ \overline u(t,\overline h(t))=0,&t\ge 0,\\[1mm]
  \overline h'(t)\geq-\mu\overline u_x(t,\overline h(t)),\ \ &t\ge 0.
 \end{array}\right.\eess
If $\overline h(0)\geq h_0$, $\overline u(0,x)\geq 0$ in $[0,\overline h(0)]$, and $\overline u(0,x)\geq u_0(x)$ in $[0,h_0]$. Then the solution $(u,h)$ of {\rm(\ref{1.1})} satisfies $h(t)\leq\overline h(t)$ in $[0,\infty)$, and $u\leq\overline u$ in $D$, where $D=\{(t,x): t\geq 0,\, 0\leq x\leq h(t)\}$.
\end{lem}

For any given $\ell>0$, let $\lambda_1(\ell;d,a)$ be the principal eigenvalue of the $T$-periodic eigenvalue problem
 \bes\left\{\begin{array}{ll}
 \phi_t-d\phi_{xx}-a(t,x)\phi=\lambda\phi, \ \ &0\le t\le T, \ \ 0<x<\ell,\\[1mm]
 B[\phi](t,0)=0, \ \ \phi(t,\ell)=0,\ \ &0\le t\le T, \\[1mm]
 \phi(0,x)=\phi(T,x),&0\le x\le\ell.
 \end{array}\right.\lbl{3.1}\ees

\begin{prop}\lbl{p3.1} \,The principal eigenvalue $\lambda_1(\ell;d,a)$ is continuous and strictly decreasing in $a$ and $\ell$. Moreover, $\lim_{\ell\to 0^+}\lambda_1(\ell;d,a)=\infty$ and
\bes
 \lim_{d\to\infty}\lambda_1(\ell;d,a)=\infty.
\lbl{3.2}\ees
\end{prop}

{\bf Proof}. We only prove (\ref{3.2}), the other conclusions can be found in the monograph \cite{Hess}. Define $\hat a=\sup_{t,x\ge 0}a(t,x)$. Then $\lambda_1(\ell;d,a)\ge\lambda_1(\ell;d,\hat a)$ since $\lambda_1(\ell;d,a)$ is decreasing in $a$. Because $\hat a$ is a constant, we know that $\lambda_1(\ell;d,\hat a)$ is the principal eigenvalue of the elliptic problem
\bess\left\{\begin{array}{ll}
 -d\phi''-\hat a\phi=\lambda\phi, \ \ 0<x<\ell,\\[1mm]
 B[\phi](0)=0, \ \ \phi(\ell)=0,
 \end{array}\right.\eess
and $\lim_{d\to\infty}\lambda_1(\ell;d,\hat a)=\infty$. Thus (\ref{3.2}) holds. This completes the proof. \ \ \ \fbox{}

In order to study the spreading phenomenon and establish the sharp criteria in later, we shall introduce some sets and analyze their properties. For any given $d>0$, define  $\sum_d=\big\{\ell>0:\,\lambda_1(\ell;d,a)=0\big\}$. By the monotonicity of $\lambda_1(\ell;d,a)$ in $\ell$, we see that $\sum_d$ contains at most one element. For any given $\ell>0$, we let $\sum_\ell^-=\big\{d>0:\,\lambda_1(\ell;d,a)\le 0\big\}$ and $\sum_\ell^+=\big\{d>0:\,\lambda_1(\ell;d,a)>0\big\}$. Here we should remark that because
$\lambda_1(\ell;d,a)$ is not monotone in $d$ (cf. Theorem 2.2 of \cite{HMP}), it is useless to define $\sum_\ell$ as the manner treating $\sum_d$.

First of all, $\sum_\ell^+\not=\emptyset$ by Proposition \ref{p3.1}.

\begin{remark}\lbl{r3.1}\, For the fixed $d>0$, because $\lim_{\ell\to 0^+}\lambda_1(\ell;d,a)=\infty$ and $\lim_{\ell\to\infty}\lambda_1(\ell;d,a):=\lambda_1(\infty; d,a)$ exists, we have that $\sum_d\not=\emptyset$ is equivalent to $\lambda_1(\infty; d,a)<0$.
\end{remark}

\begin{prop}\lbl{p3.2}\,Assume that the function $a(t,x)$ satisfies
 \vspace{-1mm}\begin{quote}
 {\bf(H1)}\, There exist $\varsigma>0$, $-2<\rho\le 0$, $k>1$ and $x_n$ satisfying $x_n\to\infty$ as $n\to\infty$, such that $a(t,x)\ge\varsigma x^\rho$ in $[0,T]\times[x_n,kx_n]$.
 \vspace{-1mm}\end{quote}
Then $\lambda_1(\infty; d,a)<0$, and so $\sum_d\not=\emptyset$ for any $d>0$.
\end{prop}

{\bf Proof}.\, Let $\lambda_1^D(kx_n;d,a)$ and $\gamma_1^D(x_n; d,a)$ be the principal eigenvalues of
 \bess\left\{\begin{array}{ll}
 \phi_t-d\phi_{xx}-a(t,x)\phi=\lambda\phi, \ \ &0\le t\le T, \ \ 0<x<kx_n,\\[1mm]
 \phi(t,0)=0, \ \ \phi(t,kx_n)=0,\ \ &0\le t\le T,\\[1mm]
 \phi(0,x)=\phi(T,x),&0\le x\le kx_n
 \end{array}\right.\eess
and
  \bess\left\{\begin{array}{ll}
 \psi_t-d\psi_{xx}-a(t,x)\psi=\gamma\psi, \ \ &0\le t\le T, \ \ x_n<x<kx_n,\\[1mm]
 \psi(t,x_n)=0, \ \ \psi(t,kx_n)=0,\ \ &0\le t\le T,\\[1mm]
 \psi(0,x)=\psi(T,x),&x_n\le x\le kx_n,
 \end{array}\right.\eess
respectively. Then $\lambda_1^D(kx_n;d,a)<\gamma_1^D(x_n; d,a)$ by Proposition \ref{p3.1}. In view of Proposition 17.7 in \cite{Hess}, $\lambda_1(kx_n;d,a)\le\lambda_1^D(kx_n;d,a)$. Thanks to $\rho\le 0$, one has $a(t,x)\ge\varsigma x^\rho\ge\varsigma k^\rho x_n^\rho$ in $[0,T]\times[x_n,kx_n]$. Since $\gamma_1^D(x_n; d,a)$ is decreasing in $a$, we have $\gamma_1^D(x_n; d,a)\le \gamma_1^D(x_n; d,\varsigma k^\rho x_n^\rho)$. Therefore,
 \bes
 \lambda_1(kx_n;d,a)<\gamma_1^D(x_n; d,\varsigma k^\rho x_n^\rho), \ \ \forall \  n\ge 1.
\lbl{3.3}\ees

Let $\psi(t,x)$ be the positive eigenfunction corresponding to $\gamma_1^D(x_n; d,\varsigma k^\rho x_n^\rho)$. Set  $y=x/x_n$ and $\Psi(t,y)=\psi(t,x)$. Then $\Psi(t,y)$ satisfies
\bess\left\{\begin{array}{ll}
 \Psi_t-dx_n^{-2}\Psi_{yy}-\varsigma k^\rho x_n^\rho\Psi=\gamma_1^D(x_n; d,\varsigma k^\rho x_n^\rho)\Psi, \ \ &0\le t\le T, \ \ 1<y<k,\\[1mm]
 \Psi(t,1)=0, \ \ \Psi(t,k)=0,\ \ &0\le t\le T,\\[1mm]
 \Psi(0,y)=\Psi(T,y),&1\le y\le k.
 \end{array}\right.\eess
Utilizing the inequality (17.6) of \cite{Hess}, we have
 \bes
 \gamma_1^D(x_n; d,\varsigma k^\rho x_n^\rho)=d\lambda^* x_n^{-2}-\varsigma k^\rho
 x_n^\rho=x_n^{-2}(d\lambda^*-\varsigma k^\rho x_n^{2+\rho})<0 \ \ \mbox{for} \ \ n\gg 1
  \lbl{3.4}\ees
since $2+\rho>0$ and $x_n\to\infty$ as $n\to\infty$, where $\lambda^*$ is the principal eigenvalue of
 \bess\left\{\begin{array}{ll}
 -u''=\lambda u, \ \ 1<y<k,\\[1mm]
 u(1)=0, \ \ u(k)=0.
 \end{array}\right.\eess
It follows from (\ref{3.3}) and (\ref{3.4}) that $\lambda_1(\infty; d,a)<0$. The proof is complete. \ \ \ \fbox{}

The condition {\bf(H1)} seems to be ``weak" because $a(t,x)$ may be ``very negative" in the sense that both $|\{a(t,x)>0\}|\ll|\{a(t,x)<0\}|$ and $\int_0^\infty a(t,x){\rm d}x=-\infty$ are allowed for any $t>0$.

\begin{prop}\lbl{p3.3}\,Assume that
 \vspace{-1mm}\begin{quote}
 {\bf(H2)}\, there exists $\hat x>0$ such that $\int_0^Ta(t,\hat x){\rm d}t>0$.
 \vspace{-1mm}\end{quote}
Then for any $\ell>\hat x$, the set $\sum_\ell^-$ is non-empty.
\end{prop}

{\bf Proof}.\, Owing to $\int_0^Ta(t,\hat x){\rm d}t>0$, there exists $0<\varepsilon<\min\{\hat x,\,\ell-\hat x\}$ such that $\int_0^T\underline a(t){\rm d}t>0$, where $\underline a(t)=\min_{x\in I}a(t,x)$ and $I:=[\hat x-\varepsilon,\hat x+\varepsilon]$. Let  $\hat\lambda$ be the principal eigenvalue of
 \bess\left\{\begin{array}{ll}
 -u''=\lambda u, \ \ \hat x-\varepsilon<x<\hat x+\varepsilon,\\[1mm]
 u(\hat x\pm\varepsilon)=0.
 \end{array}\right.\eess
Then $\hat\lambda>0$. Applying the inequality (17.6) of \cite{Hess} we have
 \[\lambda_1(\ell;d,a)\le d\hat\lambda-\frac 1T\int_0^T\underline a(t){\rm d}t.\]
Owing to $\int_0^T\underline a(t){\rm d}t>0$, there exists $d_0>0$ such that $\lambda_1(\ell;d,a)<0$ for all $0<d\le d_0$. This implies $(0,d_0)\subset\sum_\ell^-$.
\ \ \ \fbox{}

Now we consider the following $T$-periodic boundary value problem of logistic equation in a bounded interval $(0,\ell)$:
\bes\left\{\begin{array}{ll}
 v_t-d v_{xx}=a(t,x)v-b(t,x)v^2, \ \ &0< t\le T, \ 0<x<\ell,\\[1.5mm]
 B[v](t,0)=0, \ \ v(t,\ell)=\theta, &0\le t\le T,\\[1.5mm]
 v(0,x)=v(T,x), &0\le x\le \ell.
  \end{array}\right. \lbl{3.5}\ees

\begin{lem}\lbl{l3.2}\,Assume that the function $a$ satisfies {\bf(H1)} and  $\ell_0$ is the unique positive root of $\lambda_1(\ell;d,a)=0$. Then for any given $\ell>\ell_0$ and $\theta\ge \|a\|_\infty/b_1$, the problem $(\ref{3.5})$ has a unique positive solution,  where $b_1$ is given by {\bf(H)}.
\end{lem}

{\bf Proof}. The approach used in this proof can be regarded as the upper and lower solutions method. Since $\ell>\ell_0$, we have $\lambda_1(\ell;d,a)<0$. Let $\phi$ be the positive eigenfunction of (\ref{3.1}) corresponding to $\lambda_1(\ell;d,a)$. It is easy to verify that $\varepsilon\phi$ is a lower solution of (\ref{3.5}) and $\varepsilon\phi\le \theta$ provided $0<\varepsilon\ll 1$.

Let $z$ be the unique solution of the initial-boundary value problem
 \bess\left\{\begin{array}{ll}
 z_t-dz_{xx}=a(t,x)z-b(t,x)z^2, \ \ &t>0, \ 0<x<\ell,\\[1.5mm]
 B[z](t,0)=0, \ \ z(t,\ell)=\theta, &t\ge 0,\\[1.5mm]
 z(0,x)=\theta, &0\le x\le \ell.
  \end{array}\right.\eess
Then $\varepsilon\phi(t,x)\le z(t,x)\le \theta$ by the comparison principle for parabolic equations. For the integer $n\ge 0$, we define $z^n(t,x)=z(t+nT,x)$, $(t,x)\in [0,T]\times[0, \ell]$. Because $a,b$ are $T$-periodic in $t$, we see that $z^n$ satisfies
 \bess\left\{\begin{array}{lll}
 z^n_t-d z^n_{xx}=a(t,x)z^n-b(t,x)(z^n)^2,\ \ &0<t\le T,\ \ 0<x<\ell,\\[1mm]
 B[z^n](t,0)=0,\ \ z^n(t,\ell)=\theta,\ &0\le t\le T,\\[1mm]
  z^n(0,x)=z(nT,x), &0\le x\le\ell.
  \end{array}\right.\eess
Since
 \[\varepsilon\phi(0,x)=\varepsilon\phi(T,x)\le z(T,x)=z^1(0,x),\ \ \ z^1(0,x)=z(T,x)\le \theta=z(0,x),\] we can apply the comparison principle to get that $\varepsilon\phi\le z^1\leq z$ in $[0,T]\times[0, \ell]$. Which also implies
  $$\varepsilon\phi(0,x)\le z^2(0,x)=z^1(T,x)\le z(T,x)=z^1(0,x).$$
As above, $\varepsilon\phi\le z^2\leq z^1$ in $[0,T]\times[0, \ell]$.
Applying the inductive method we have that $z^n$ is decreasing in $n$ and $z^n\ge\varepsilon\phi$ in $[0,T]\times[0, \ell]$.
So, there exists a non-negative function $\overline v$ such that $z^n\to\overline v$ pointwise in $[0,T]\times[0,\ell]$ as $n\to\infty$.  Since $z^{n+1}(0,x)=z^n(T,x)$, it follows that $\overline v(0,x)=\overline v(T,x)$. Based on the regularity theory for parabolic equations and compact argument, it can be proved that there exists a subsequence $\{n_i\}$, such that $z^{n_i}\to\overline v$ in $C^{1,2}([0,T]\times[0,\ell])$ as $i\to\infty$, $\overline v\ge\varepsilon\phi$ in $[0,T]\times[0, \ell]$ and $\overline v$ satisfies the first two equations of (\ref{3.5}). This shows that $\overline v$ is a positive solution of (\ref{3.5}).

Now we prove that $\overline v$ is the unique positive solution of (\ref{3.5}). Let $v$ be another one. By the maximum principle, we have $v(t,x)\le \theta$, thereby, $v\le z$. This implies $v\le z^n$ in $[0,T]\times[0, \ell]$ for any integer $n\ge 0$. Certainly, $v\le \overline v$ in $[0,T]\times[0, \ell]$. This suggests that $\overline v$ is a maximal positive solution of (\ref{3.5}). It is easy to see that there exists a constant $\varepsilon>0$ such that $v\ge\varepsilon \overline v$ in $[0,T]\times[0,\ell]$, and hence the infimum
  \[\sigma=\inf_{0\le t\le T\atop 0<x<\ell}\frac{v(t,x)}{\overline v(t,x)}\]
exists and is positive. Clearly, $\sigma\le 1$ and $v\geq\sigma \overline v$ in $[0,T]\times[0,\ell]$. If we can show $\sigma=1$ then $v=\overline v$ and the uniqueness is derived. Assume on the contrary that $\sigma<1$. Denote $\varphi=v-\sigma \overline v$, then $\varphi\ge 0$, and
 \bess
&\varphi(0,x)=\varphi(T,x) \ \ \mbox{for} \ \ 0\le x\le \ell,&\\
 &B[\varphi](t,0)=0, \ \ \varphi(t,\ell)=(1-\sigma)\theta>0 \ \
 \mbox{for} \ \ 0\le t\le T.&\eess
The direct calculation yields
 \bess
 \varphi_t-d \varphi_{xx}&=&a(t,x)\varphi-b(t,x)\left(v^2-\sigma \overline v^2\right)\\[1mm]
 &>&a(t,x)\varphi-b(t,x)\left(v^2-\sigma^2 \overline v^2\right)\\[1mm]
 &\ge&a(t,x)\varphi-2b(t,x)v \varphi.\eess
The maximum principle allows us to deduce that $\varphi>0$ in $[0,T]\times(0,\ell]$. Similarly to the above, there exists $\varepsilon_1>0$ such that $\varphi\ge\varepsilon_1 \overline v$, and thus $v\ge(\sigma+\varepsilon_1)\overline v$ in $[0,T]\times[0,\ell]$. This is a contradiction with the definition of $\sigma$, and the uniqueness is derived. \ \ \ \fbox{}

\section{Spreading-vanishing dichotomy and long time behavior of solution}
\setcounter{equation}{0}

We first give a lemma, by which the vanishing phenomenon is immediately obtained. Moreover, this lemma will play an important role in the establishment of sharp criteria for spreading and vanishing.

\begin{lem}\lbl{l4.1}\, Let $d,\mu$ and $B$ be as above, $c\in\mathbb{R}$. Assume that $g\in C^1([0,\infty))$, $\varphi\in C^{\frac{1+\nu}2,1+\nu}([0,\infty)\times[0,g(t)])$ and satisfy $g(t)>0$, $\varphi(t,x)>0$ for $t\ge 0$ and $0<x<g(t)$. We further suppose that
$\lim_{t\to\infty} g(t)<\infty$, $\lim_{t\to\infty} g'(t)=0$ and there exists a constant $C>0$ such that $\|\varphi(t,\cdot)\|_{C^1[0,\,g(t)]}\leq C$ for $t>1$. If $(\varphi,g)$ satisfies
  \bess\left\{\begin{array}{lll}
 \varphi_t-d \varphi_{xx}\geq c \varphi, &t>0,\ \ 0<x<g(t),\\[.5mm]
 B[\varphi]=0, \ &t\ge 0, \ \ x=0,\\[.5mm]
 \varphi=0,\ \ g'(t)\geq-\mu \varphi_x, \ &t\ge 0,\ \ x=g(t),
 \end{array}\right.\eess
then $\lim_{t\to\infty}\max_{0\leq x\leq g(t)}\varphi(t,x)=0$.
 \end{lem}

{\bf Proof}.\, When $\alpha=0$ or $\beta=0$, this is exactly Proposition 3.1 of \cite{Wjde}. When $\alpha>0$ and $\beta>0$, that proof remains valid. The details are omitted here. \ \ \ \fbox{}

If $h_\infty<\infty$, the estimate (\ref{2.3}) implies $\lim_{t\to\infty} h'(t)=0$. Applying (\ref{2.3}), (\ref{2.4}) and Lemma \ref{l4.1}, we have the following theorem.

\begin{theo}\lbl{th4.1}{\rm(}Vanishing{\rm)}\, Let $(u,h)$ be the solution of {\rm(\ref{1.1})}. When $h_\infty<\infty$, we must have
 \bes
 \lim_{t\to\infty}\,\max_{0\leq x\leq h(t)}u(t,x)=0.\lbl{4.1}
 \ees
This shows that if the species cannot spread successfully, it will be extinct in the long run.
\end{theo}

In the following we investigate the spreading phenomenon. To this aim, we first study the existence and uniqueness of positive solution to the following $T$-periodic problem
 \bes\left\{\begin{array}{ll}
 U_t-dU_{xx}=a(t,x)U-b(t,x)U^2, \ \ &0\le t\le T, \ 0<x<\infty,\\[1mm]
 B[U](t,0)=0,  &0\le t\le T,\\[1mm]
 U(0,x)=U(T,x), &0\le x<\infty.
\end{array}\right. \lbl{4.2}\ees

\begin{theo}\lbl{th4.2}\, Assume that there exist a constant $\rho$ with $-2<\rho\le 0$, and $T$-periodic positive functions $a_\infty(t),b_\infty(t),a^\infty(t),b^\infty(t)\in C^{\nu/2}([0,T])$, such that
 \bes\left\{\begin{array}{ll}
 a_\infty(t)=\displaystyle\liminf_{x\to\infty}\frac{a(t,x)}{x^\rho}, \ \ \
 a^\infty(t)=\limsup_{x\to\infty}\frac{a(t,x)}{x^\rho},\\[2mm]
 b_\infty(t)=\displaystyle\liminf_{x\to\infty}b(t,x), \ \ \
 b^\infty(t)=\limsup_{x\to\infty}b(t,x)
 \end{array}\right.\lbl{4.3}\ees
uniformly on $[0,T]$. Then $(\ref{4.2})$ has a unique positive solution $U\in C^{1+\frac{\nu}2, 2+\nu}([0,T]\times[0,\infty))$, and satisfies
 \bes
 \frac{\min_{[0,T]}a_\infty(t)}{\max_{[0,T]}b^\infty(t)}\le
 \liminf_{x\to\infty}\frac{U(t,x)}{x^\rho}, \ \ \ \limsup_{x\to\infty}\frac{U(t,x)}{x^\rho}\le \frac{\max_{[0,T]}a^\infty(t)}{\min_{[0,T]}b_\infty(t)}
\lbl{4.4}\ees
uniformly on $[0,T]$.

It is worth mentioning that the condition {\rm(\ref{4.3})} implies the assumptions {\bf(H1)} and {\bf(H2)}.
 \end{theo}

{\bf Proof}.\, When $\alpha=0$, this theorem is a special case of Theorem 1.3 in \cite{PD}. We consider the case $\alpha>0$ in the following. This proof is divided into three steps. In the first one, we construct the minimal positive solution of (\ref{4.2}). The estimate (\ref{4.4}) will be given in the second step. Finally, we show the uniqueness of positive solution.

{\it Step 1: The existence}. In this step, we shall construct a positive solution $\underline U$ and prove that it is the minimal one. Consider the following auxiliary problem
\bes\left\{\begin{array}{ll}
 U_t-d U_{xx}=a(t,x)U-b(t,x)U^2,\ \ &0\le t\le T, \ \ 0<x<\ell,\\[1mm]
  B[U](t,0)=U(t,\ell)=0, &0\le t\le T,\\[1mm]
 U(0,x)=U(T,x), &0\le x\le \ell.\end{array}\right.
 \lbl{4.5}\ees
Let $\lambda_1(\ell;d,a)$ be the principal eigenvalue of (\ref{3.1}). Since the assumption {\bf(H1)} holds, by use of Proposition \ref{p3.2}, there exists $\ell_0\gg 1$ such that $\lambda_1(\ell,d,a)<0$ for all $\ell\ge\ell_0$. For such $\ell$, utilizing Theorem 28.1 of \cite{Hess}, the problem (\ref{4.5}) admits a unique positive solution, denoted by $U_\ell(t,x)$. Obviously, $U_\ell\le\|a\|_\infty/b_1$ by the maximum principle, where $b_1$ is given by {\bf(H)}. For $\ell^*>\ell$, it is easy to see that $U_{\ell^*}$ is an upper solution of (\ref{4.5}). Let $\phi(t,x)$ be the positive eigenfunction of (\ref{3.1}) corresponding to $\lambda_1(\ell,d,a)$ and  $\varepsilon>0$ be a constant. Then $\varepsilon\phi$ is a positive lower solution of (\ref{4.5}) and $\varepsilon\phi\le U_{\ell^*}$ provided $\varepsilon\ll 1$. Thus, $U_{\ell^*}\geq U_\ell$ since $U_\ell$ is the unique positive solution of (\ref{4.5}). This shows that $U_\ell$ is increasing in $\ell$. Make use of the regularity theory for parabolic equations and compact argument, it can be proved that there exists a subsequence of $\{U_\ell\}$, denoted by itself, and a positive function $\underline U\in C^{1, 2}([0,T]\times[0,\infty))$, such that $U_\ell\to\underline U$ in $C^{1, 2}([0,T]\times[0,L])$ for any $L>0$, and $\underline U$ solves (\ref{4.2}).

Let $U$ be a positive solution of (\ref{4.2}). Then $U\le\|a\|_\infty/b_1$ by the maximum principle. Obviously, $U$ is an upper solution of (\ref{4.5}) for any given $\ell\ge\ell_0$. As above, $U_\ell\le U$ in $[0,T]\times[0,\ell]$. This implies $\underline U\le U\le\|a\|_\infty/b_1$ in $[0,T]\times[0,\infty)$.

{\it Step 2: Proof of {\rm(\ref{4.4})}}.\, For any positive solution $U$ of (\ref{4.2}), we have known that $\underline U\le U\le\|a\|_\infty/b_1$ in $[0,T]\times[0,\infty)$. Take account of Theorem 1.3 in \cite{PD}, the problem
  \bes\left\{\begin{array}{ll}
 W_t-dW_{xx}=a(t,x)W-b(t,x)W^2, \ \ &0\le t\le T, \ 0<x<\infty,\\[1mm]
 W_x(t,0)=0,  &0\le t\le T,\\[1mm]
 W(0,x)=W(T,x), &0\le x<\infty
\end{array}\right. \lbl{4.6}\ees
has a unique positive solution $W(t,x)$. For $\ell>0$, let us consider the problem
 \bes\left\{\begin{array}{ll}
 V_t-d V_{xx}=a(t,x)V-b(t,x)V^2, \ \ &0\le t\le T, \ 0<x<\ell,\\[1mm]
 V_x(t,0)=0, \ \ V(t,\ell)=U(t,\ell), \ \  &0\le t\le T,\\[1mm]
 V(0,x)=V(T,x), &0\le x\le\ell.
 \end{array}\right. \lbl{4.7}\ees
It is obvious that $U$ is a lower solution of (\ref{4.7}) since $U_x(t,0)>0$, and the constant $K=1+\|a\|_\infty/b_1$ is an upper solution of (\ref{4.7}). Similarly to the proof of Lemma \ref{l3.2}, we can prove that the problem (\ref{4.7}) has a unique  positive solution $V_\ell$ and $U\le V_\ell\le K$ in $[0,T]\times[0,\ell]$.

Arguing as Step 1, there exists a subsequence of $\{V_\ell\}$, denoted by itself, such that $V_\ell\to W$ in $C^{1, 2}([0,T]\times[0,L])$ for any $L>0$ since $W$ is the unique positive solution of (\ref{4.6}). Therefore, $U\le W$ in $[0,T]\times[0,\infty)$. In view of Eq. (2.6) in \cite{PD}, we have that
 \bes
 \limsup_{x\to\infty}\frac{U(t,x)}{x^\rho}\le\limsup_{x\to\infty}\frac{W(t,x)}{x^\rho}
 \le \frac{\max_{[0,T]}a^\infty(t)}{\min_{[0,T]}b_\infty(t)}\lbl{4.8}\ees
uniformly on $[0,T]$.

Set
 \[\underline a(x)=\min_{[0,T]}a(t,x), \ \ \overline b(x)=\max_{[0,T]}b(t,x),\]
and consider the problem
 \bes\left\{\begin{array}{ll}
 -dw''=\underline a(x)w-\overline b(x)w^2, \ \ &0<x<\infty,\\[1mm]
 w(0)=0.
 \end{array}\right.\lbl{4.9}\ees
Thanks to the conditions {\bf(H)} and (\ref{4.3}), we can show that $\underline a(x),\overline b(x)\in C^\nu([0,\infty))$, and
 \bess\left\{\begin{array}{ll}
 0<\displaystyle\min_{[0,T]}a_\infty(t)=\liminf_{x\to\infty}\frac{\underline a(x)}{x^\rho}\le
 \limsup_{x\to\infty}\frac{\underline a(x)}{x^\rho}\le \max_{[0,T]}a^\infty(t),\\[3mm]
 0<\displaystyle\min_{[0,T]}b_\infty(t)\le\liminf_{x\to\infty}\overline b(x)\le\limsup_{x\to\infty}\overline b(x)
=\max_{[0,T]}b^\infty(t).
 \end{array}\right.\eess
By virtue of Proposition \ref{p3.2}, the first eigenvalue of
 \bess\left\{\begin{array}{ll}
 -d\phi''-\underline a(x)\phi=\lambda\phi, \ \ 0<x<\ell,\\[1mm]
 \phi(0)=0, \ \ \phi(\ell)=0
 \end{array}\right.\eess
is negative provided $\ell\gg 1$. And then, the problem
\bess\left\{\begin{array}{ll}
 -dw''=\underline a(x)w-\overline b(x)w^2, \ \ &0<x<\ell,\\[1mm]
 w(0)=0, \ \ w(\ell)=0
 \end{array}\right. \eess
has a unique positive solution $w_\ell(x)$. By the same argument as in Step 1, we can show that there exists a positive function $w$ such that $w_\ell\to w$ in $C^{1, 2}([0,L])$ for any $L>0$, and $w$ solves (\ref{4.9}). Obviously, $w_\ell$ is a lower solution of (\ref{4.5}). Hence, $w_\ell\leq U_\ell$ since $U_\ell$ is the unique positive solution of (\ref{4.5}). This implies $w\le\underline U \le U$ in $[0,T]\times[0,\infty)$ since $w_\ell\to w$ and $U_\ell\to\underline U$. Moreover, $w'(0)>0$ by the uniqueness.

If $w'(x)>0$ in $[0,\infty)$, then $w(x)\to w_*$ as $x\to\infty$ for some positive constant $w_*$. As $\rho\le 0$, it is immediate to get
 \bes
 \liminf_{x\to\infty}\frac{U(t,x)}{x^\rho}\ge\liminf_{x\to\infty}\frac{w(x)}{x^\rho}\ge
 w_*\lbl{4.10}\ees
uniformly on $[0,T]$. In this case we claim that $\rho=0$ and
 \bes
 w_*\geq\frac{\liminf_{x\to\infty}\underline a(x)}{\limsup_{x\to\infty}\overline b(x)}
 =\frac{\min_{[0,T]}a_\infty(t)}{\max_{[0,T]}b^\infty(t)}.\lbl{4.11}\ees
In fact, if $\rho<0$, then $\lim_{x\to\infty}a(t,x)=0$ uniformly on $[0,T]$ by (\ref{4.3}). Therefore, $\lim_{x\to\infty}\underline a(x)=0$. This is impossible since $\overline b(x)$ has a positive lower bound, $w'(x)>0$ in $[0,\infty)$ and $w$ is bounded from above. So, $\rho=0$. On the contrary we assume that (\ref{4.11}) does not hold. Denote $a_*=\liminf_{x\to\infty}\underline a(x)$ and $b^*=\limsup_{x\to\infty}\overline b(x)$. Then there exist $\varepsilon>0$ and $x^*\gg 1$, such that $w(x)<{a_*}/{b^*}-\varepsilon$ for all $x\ge x^*$.
For such $\varepsilon>0$, there exist $\delta>0$ and $x_0>x^*$ such that
 \bess
 \frac{a_*-\delta}{b^*+\delta}-\frac{a_*}{b^*}+\varepsilon:=\sigma>0,
 \eess
and $\underline a(x)>a_*-\delta,\,\overline b(x)<b^*+\delta$ for all $x>x_0$. Then
 \bess
-dw''=w\big[\underline a(x)-\overline b(x)w\big]>w\overline b(x)\left(\frac{\underline a(x)}{\overline b(x)}-\frac{a_*}{b^*}+\varepsilon\right)>\sigma\overline b(x)w, \ \ \forall \  x>x_0,
\eess
this is impossible since $\sigma>0$ and $\overline b(x)$ has a positive lower bound. So,  (\ref{4.11}) is true.

If $w'(x_0)=0$ for some $x_0>0$, we set
 \bess
 &\underline a^*(x)=\underline a(x+x_0), \ \ \overline b^*(x)=\overline b(x+x_0), \ \ v(x)= w(x+x_0) \ \ \mbox{for} \ \ x\geq 0,&\\[1mm]
 &\underline a^*(x)=\underline a(-x+x_0), \ \  \overline b^*(x)=\overline b(-x+x_0), \ \ v(x)=w(-x+x_0) \ \
  \mbox{for}\ \ x< 0.\eess
Then the function $v(x)$ satisfies
 \bes
 -dv''=\underline a^*(x)v-\overline b^*(x)v^2, \ \ -\infty<x<\infty.
 \lbl{4.12}\ees
Applying Theorem 7.12 of \cite{Du}, we have that $v(x)$ is the unique positive solution of (\ref{4.12}) and satisfies
\bess
\liminf_{x\to\infty}\frac{v(x)}{x^\rho}\ge\frac{\displaystyle\liminf_{x\to\infty}\frac{\underline a^*(x)}{x^\rho}}{\displaystyle\limsup_{x\to\infty}\overline b^*(x)}=\frac{\displaystyle\liminf_{x\to\infty}\frac{\underline a(x)}{x^\rho}}{\displaystyle\limsup_{x\to\infty}\overline b(x)}=\frac{\min_{[0,T]}a_\infty(t)}{\max_{[0,T]}b^\infty(t)}.
 \eess
Therefore,
 \bes
\liminf_{x\to\infty}\frac{U(t,x)}{x^\rho}\ge\liminf_{x\to\infty}\frac{w(x)}{x^\rho}=\liminf_{x\to\infty}\frac{v(x)}{x^\rho}
\ge\frac{\min_{[0,T]}a_\infty(t)}{\max_{[0,T]}b^\infty(t)}
 \lbl{4.13}\ees
uniformly on $[0,T]$.

It follows from (\ref{4.8}), (\ref{4.10}), (\ref{4.11}) and (\ref{4.13}) that $U$ satisfies (\ref{4.4}).

{\it Step 3: The uniqueness}. Let $U$ be a positive solution of (\ref{4.2}). Then  $U\ge\underline U$, here $\underline U$ is the minimal positive solution of (\ref{4.2}) obtained by Step 1. To prove the uniqueness, it suffices to show that $U\equiv\underline U$. If this is not true, then $U\ge,\,\not=\underline U$. It follows from (\ref{4.4}) that there exists $k>1$ such that
$U\le k\underline U$ in $[0,T]\times[0,\infty)$. To arrive at a contradiction, we turn to a technique introduced by Marcus and V\'{e}ron in \cite{MV}. Define
$V=\underline U-(2k)^{-1}(U-\underline U)$. Then
 \bes
 \underline U>V\ge\frac{k+1}{2k}\underline U, \ \ \ \frac{2k}{2k+1}V+\frac 1{2k+1}U=\underline U.
 \lbl{4.14}\ees
Noticing that $z^2$ is convex in $z\in(0,\infty)$, we have
$\underline U^2\le \frac{2k}{2k+1}V^2+\frac 1{2k+1}U^2$ by the second formula of (\ref{4.14}).
Because $b(t,x)$ is positive, then the direct computation gives
 \[V_t-dV_{xx}\ge a(t,x)V-b(t,x)V^2.\]
Obviously, for the large $\ell$,
 \[B[V](t,0)=0, \ \ V(t,\ell)>0 \ \ \mbox{in}\ \ [0, T]; \ \ \
 V(0,x)=V(T,x) \ \ \mbox{in} \ \ [0,\ell]. \]
This shows that $V$ is an upper solution of (\ref{4.5}). Let $U_\ell$ be the unique positive solution of (\ref{4.5}). As above, $U_\ell\le V$ in $[0,T]\times[0,\ell]$. Hence, $\underline U\leq V$ in $[0,T]\times[0,\infty)$ due to $U_\ell\to\underline U$ as $\ell\to\infty$. This is a contradiction with the first inequality of (\ref{4.14}). So, $U\equiv\underline U$ and the uniqueness is derived. The proof is complete. \ \ \ \fbox{}

\begin{theo}\lbl{th4.3}{\rm(}Spreading{\rm)}\, Assume that $(\ref{4.3})$ holds. If $h_\infty=\infty$, then $\lim_{n\to\infty}u(t+nT,x)=U(t,x)$ uniformly in $[0,T]\times[0,L]$ for any $L>0$.
\end{theo}

{\bf Proof}.\, At present, the condition {\bf(H1)} holds. Let $\lambda_1(\ell;d,a)$ be the principal eigenvalue of (\ref{3.1}). In view of Proposition \ref{p3.2}, we can choose an $\ell_0\gg 1$ such that $\lambda_1(\ell,d,a)<0$ for all $\ell\ge\ell_0$. We divide the proof into two parts, and show that, respectively,
\bes
 \limsup_{n\to\infty}u(t+nT,x)\leq U(t,x) \ \ \ \mbox{uniformly in} \ \ [0,T]\times[0,L]
\lbl{4.15}\ees
and
 \bes
 \liminf_{n\to\infty}u(t+nT,x)\geq U(t,x) \ \ \ \mbox{uniformly in} \ \ [0,T]\times[0,L].\lbl{4.16}\ees

{\it Step 1}. Take $\theta=\|u_0\|_\infty+\frac{1}{b_1}\|a\|_\infty$, where $b_1$ is given by {\bf(H)}. Then $u<\theta$ by the maximum principle. For the fixed $\ell>\ell_0$, there exists an integer $m\gg 1$ such that $h(t)>\ell$ for all $t\ge mT$.

Since $\lambda_1(\ell;d,\|a\|_\infty)\le\lambda_1(\ell;d,a)<0$, by the upper and lower solutions method we can prove that the boundary value problem
  \bess\left\{\begin{array}{ll}
 -d V_{xx}=\|a\|_\infty V-b_1 V^2, \ \ &0<x<\ell,\\[1.5mm]
 B[V](0)=0, \ \ V(\ell)=\theta
  \end{array}\right.\eess
has a unique positive solution $V(x)$. In consideration of the regularity of $u(mT,x)$ and $V(x)$ in $x$, we can find a constant $k\ge 1$ such that $u(mT,x)\le k V(x)$ for all $0\le x\le\ell$. Note that $u(t,\ell)<k\theta=kV(\ell)$, we can apply the comparison principle to $u$ and $kV$, and then derive that $u(t,x)\le kV(x)$ for all $t\ge mT$ and $0\le x\le\ell$. Since $k\ge 1$, it is easy to see that the function $v:=kV$ satisfies
  \[-d v_{xx}\ge\|a\|_\infty v-b_1 v^2\ge a(t,x)v-b(t,x)v^2.\]

Let $w_\ell$ be the unique solution of
 \bess\left\{\begin{array}{lll}
 w_t-dw_{xx}=a(t,x)w-b(t,x)w^2,\ \ &t>mT,\ \ 0<x<\ell,\\[1mm]
 B[w](t,0)=0,\ \ w(t,\ell)=k V(\ell),\ &t>mT,\\[1mm]
  w(mT,x)=k V(x), &0\le x\le\ell.
  \end{array}\right.\eess
The comparison principle gives $u\leq w_\ell\le kV$ for $t\ge mT$ and $0\leq x\leq \ell$. For the integer $n\ge m$, we define $w_\ell^n(t,x)=w_\ell(t+nT,x)$, $(t,x)\in [0,T]\times[0, \ell]$. Similarly to the proof of Lemma \ref{l3.2}, it can be shown that $w^{n}_\ell\to W_\ell$ in $C^{1,2}([0,T]\times[0, \ell])$ as $n\to\infty$,
where $W_\ell$ is the unique positive solution of (\ref{3.5}) with $\theta$ replaced by $k\theta$. Owing to $u(t+nT,x)\le w_\ell(t+nT,x)=w^{n}_\ell(t,x)$ in $[0,T]\times[0, \ell]$, we get
 \bes
 \limsup_{n\to\infty}u(t+nT,x)\leq W_\ell(t,x) \ \ \ \mbox{uniformly on} \ \ [0,T]\times[0,\ell].
\lbl{4.17}\ees

Since $W_\ell\le k\theta$, it is easy to see that $W_\ell$ is decreasing in $\ell$. Remember that $U$ is the unique positive solution of (\ref{4.2}), arguing as Step 1 in the proof of Theorem \ref{th4.2}, we can prove that $W_\ell\to U$ in $C^{1, 2}([0,T]\times[0,L])$ as $\ell\to\infty$. Combining this fact with (\ref{4.17}), we get (\ref{4.15}).

{\it Step 2}. Let $\ell>\ell_0$ and $v_\ell$ be the unique solution of
 \bes\left\{\begin{array}{lll}
 v_t-d v_{xx}=a(t,x)v-b(t,x)v^2,\ \ &t>mT, \ \ 0<x<\ell,\\[1mm]
 B[v](t,0)=0, \ \ v(t,\ell)=0,\ \ \ &t>mT,\\[1mm]
 v(mT,x)=u(mT,x),&x\in [0,\ell].
 \end{array}\right. \lbl{4.18}\ees
Then $u\geq v_\ell$ in $[mT,\infty)\times[0,\ell]$. Since $\lambda_1(\ell;d,a)<0$, the problem
 \bess\left\{\begin{array}{ll}
 V_t-d V_{xx}=a(t,x)V-b(t,x)V^2, \ \ &0\le t\le T, \ 0<x<\ell,\\[1mm]
 B[V](t,0)=0, \ \ V(t,\ell)=0, \ \  &0\le t\le T,\\[1mm]
 V(0,x)=V(T,x), &0\le x\le\ell
 \end{array}\right.\eess
admits a uniqe positive solution $V_\ell$. Clearly, the constant $\theta=\|u_0\|_\infty+\frac{1}{b_1}\|a\|_\infty$ is an upper solution of (\ref{4.18}). Noticing  $v_\ell(mT,x)=u(mT,x)>0$ in $(0,\ell]$, take advantage of Theorem 28.1 in \cite{Hess}, it follows that $v_\ell(t+nT,x)\to V_\ell(t,x)$ in $C^{1,2}([0,T]\times[0, \ell])$ as $n\to\infty$.

It is obvious that $V_\ell$ is increasing in $\ell$. Similarly to Step 1, we can derive that $V_\ell\to U$ in $C^{1, 2}([0,T]\times[0,L])$. Consequently, the limit (\ref{4.16}) is obtained since $u\geq v_\ell$ in $[mT,\infty)\times[0,\ell]$ for all $\ell>\ell_0$. The proof is finished.\ \ \ \fbox{}

\section{Sharp criteria for spreading and vanishing, spreading speed}
\setcounter{equation}{0}{\setlength\arraycolsep{2pt}

We first give a necessary condition for vanishing. Let $\lambda_1(\ell;d,a)$ be the principal eigenvalue of (\ref{3.1}).

\begin{lem}\lbl{l5.1}\, If $h_\infty<\infty$, then $\lambda_1(h_\infty;d,a)\geq 0$.
\end{lem}

{\bf Proof}.\, We assume $\lambda_1(h_\infty;d,a)<0$ to get a contradiction. By the continuity of $\lambda_1(\ell;d,a)$ in $\ell$ and $h(t)\to h_\infty$, there exists $\tau\gg 1$ such that $\lambda_1(h(\tau);d,a)<0$. Let $w$ be the unique solution of
  $$\left\{\begin{array}{ll}
 w_t-dw_{xx}=a(t,x)w-b(t,x)w^2, \ \ &t>\tau, \ \, 0<x<h(\tau),\\[1mm]
 B[w](t,0)=w(t,h(\tau))=0,&t\ge\tau,\\[1mm]
  w(\tau,x)=u(\tau,x),&0\le x\le h(\tau).
  \end{array}\right.$$
Then $u\geq w$ in $[\tau,\infty)\times[0,h(\tau)]$ by the comparison principle. Remembering $\lambda_1(h(\tau);d,a)<0$, it follows from Theorem 28.1 of \cite{Hess} that $w(t+nT,x)\to z(t,x)$ uniformly in $[0,T]\times[0,h(\tau)]$ as $n\to\infty$, where $z$ is the unique positive solution of the following $T$-periodic boundary value problem
\bess\left\{\begin{array}{ll}
 z_t-dz_{xx}=a(t,x)z-b(t,x)z^2,\ \ &0\le t\le T, \ 0<x<h(\tau),\\[1mm]
  B[z](t,0)=z(t,h(\tau))=0, &0\le t\le T,\\[1mm]
 z(0,x)=z(T,x), &0\le x\le h(\tau).\end{array}\right.
 \eess
Since $u\geq w$ in $[\tau,\infty)\times[0,h(\tau)]$, it is deduced immediately that
 \[\liminf_{n\to\infty}u(t+nT,x)\ge z(t,x), \ \ \forall \  (t,x)\in [0,T]\times[0,h(\tau)].\]
This is a contradiction with (\ref{4.1}). The proof is complete. \ \ \ \fbox{}

\begin{lem}\label{l5.2}\, If $\lambda_1(h_0;d,a)>0$, then there exists $\mu_0>0$, such that $h_\infty<\infty$ provided $\mu\leq\mu_0$. Hence, by Lemma {\rm\ref{l5.1}}, $\lambda_1(h_\infty;d,a)\geq 0$ for $\mu\leq\mu_0$.
\end{lem}

{\bf Proof}.\,Let $\phi(t,x)$ be the corresponding positive eigenfunction to $\lambda_1:=\lambda_1(h_0;d,a)$ of (\ref{3.1}) with $\ell=h_0$. Noticing that $\phi_x(t,h_0)<0$, $\phi(t,0)>0$ in $[0,T]$ when $\beta>0$, and $\phi_x(t,0)>0$ in $[0,T]$ when $\beta=0$. By the regularity of $\phi$, there exists a constant $C>0$ such that
 \bes
 x\phi_x(t,x)\le C\phi(t,x) , \ \ \forall \ (t,x)\in[0,T]\times[0, h_0].
 \lbl{5.1}\ees
Let $0<\delta,\,\sigma<1$ and $K>0$ be constants, which will be determined later. Set
 \bess
 \displaystyle s(t)&=&1+2\delta-\delta {\rm e}^{-\sigma t}, \ \ \tau(t)=\int_0^t s^{-2}(\rho){\rm d}\rho, \ \ t\geq 0, \\[.1mm]
 v(t,x)&=&K{\rm e}^{-\sigma t}\phi\left(\tau(t),\,y\right), \ \ y=y(t,x)=\frac x{s(t)}, \ \ 0\leq x\leq h_0s(t).
 \eess

Firstly, for any given $0<\varepsilon\ll 1$, since $a$ is uniformly continuous in $[0,T]\times[0,3h_0]$ and $T$-periodic in $t$, we have that there exists $0<\delta_0(\varepsilon)\ll 1$ such that, for all $0<\delta\le\delta_0(\varepsilon)$ and $0<\sigma<1$,
 \bes
 \left|s^{-2}(t)a\big(\tau(t),\,y(t,x)\big)-a(t,x)\right|\le\varepsilon, \ \ \forall \ t>0, \ \ 0\leq x\leq h_0s(t).
  \lbl{5.2}\ees
Note that (\ref{5.1}), (\ref{5.2}) and $\lambda_1>0$, the direct calculation yields,
 \bes
 v_t-dv_{xx}-a(t,x)v&=&v\left(-\sigma+\frac {a(\tau,y)}{s^2(t)}-a(t,x)
-\frac{y\phi_y(\tau,y)}{\phi(\tau,y)}\frac{\sigma\delta}{s(t)}{\rm e}^{-\sigma t}+\frac{\lambda_1}{s^2(t)}\right)\nonumber\\[2mm]
 &\geq&v(-\sigma-\varepsilon-C\sigma+\lambda_1/4)
 >0, \ \ \forall \ t>0, \ \ 0<x<h_0s(t)\qquad
 \lbl{5.3}\ees
provided $0<\sigma,\varepsilon\ll 1$.

\vskip 2pt Evidently, $v(t,h_0s(t))=K{\rm e}^{-\sigma t}\phi(\tau(t),h_0)=0$. If either $\alpha=0$ or $\beta=0$,  then $B[v](t,0)=0$. If $\alpha,\beta>0$, then $\alpha\phi(\tau(t),0)=\beta\phi_y(\tau(t),0)$, and $\phi_y(\tau(t),0)>0$ by the Hopf lemma. Therefore,
 \[B[v](t,0)=\beta K{\rm e}^{-\sigma t}\phi_y(\tau(t),0)[1-1/s(t)]>0\]
due to $s(t)>1$. In a word,
 \bes
 B[v](t,0)\ge 0, \ \ v(t,h_0s(t))=0, \ \ \forall \  t>0.
 \lbl{5.4}\ees
Fix $0<\sigma,\varepsilon\ll 1$ and $0<\delta\le\delta_0(\varepsilon)$. Based on the regularities of $u_0(x)$ and $\phi(x)$, we can choose a constant $K\gg1$ such that
 \bes
 u_0(x)\leq K\phi\left(0,\,x/(1+\delta)\right)=v(0,x), \ \ \forall \ 0\le x\le h_0.
\lbl{5.5}\ees
Thanks to $h_0s'(t)=h_0\sigma\delta{\rm e}^{-\sigma t}$ and $v_x(t,h_0s(t))=\frac 1{s(t)}K{\rm e}^{-\sigma t}\phi_y(\tau(t),h_0)$, there exists $\mu_0>0$ such that, for all $0<\mu\le \mu_0$,
 \bes
 h_0s'(t)\geq-\mu v_x(t,h_0s(t)), \ \ \forall \ t\ge 0.
\lbl{5.6}\ees

Remember (\ref{5.3})-(\ref{5.6}). Applying Lemma \ref{l3.1} to $(u,h(t))$ and $(v,h_0s(t))$, it follows that
 \[h(t)\leq h_0s(t), \ \ u(t,x)\leq v(t,x), \ \ \forall \  t\geq 0, \ \ 0\leq x\leq h(t).\]
Hence $h_\infty\leq h_0s(\infty)=h_0(1+2\delta)$ for all $0<\mu\leq\mu_0$. This finishes the proof. \ \ \ \fbox{}

\begin{lem}\lbl{l5.3}\, Let $C>0$ be a constant. For any given constants $\overline h_0, H>0$, and any function $\overline u_0\in C^2([0,\overline h_0])$ satisfying $B[\overline u_0](0)=\overline u_0(\overline h_0)=0$ and $\overline u_0>0$ in $(0,\overline h_0)$, there exists $\mu^0>0$ such that when $\mu\geq\mu^0$ and $(\overline u, \overline h)$ satisfies
 \bess
 \left\{\begin{array}{ll}
   \overline u_t-d\overline u_{xx}\geq -C \overline u, \ &t>0, \ \ 0<x< \overline h(t),\\[1mm]
  B[\overline u](t,0)=0=\overline u(t, \overline h(t)),\ &t\geq 0,\\[1mm]
 \overline h'(t)=-\mu \overline u_x(t, \overline h(t)), \ &t\geq 0,\\[1mm]
  \overline u(0,x)=\overline u_0(x),\ \ \overline h(0)=\overline h_0, \ &0\leq x\leq \overline h_0,
 \end{array}\right.
 \eess
we must have $\liminf_{t\to\infty}\overline h(t)>H$.
\end{lem}

The proof of Lemma \ref{l5.3} is essentially similar to that of Lemma 3.2 in \cite{WZjdde} and is hence omitted.

Now we fix $d$, and consider $h_0$ and $\mu$ as varying parameters to depict the sharp criteria for spreading and vanishing. Assume that $\sum_d\not=\emptyset$ and let $h^*=h^*(d)\in \sum_d$, i.e., $\lambda_1(h^*;d,a)=0$. It is worth mentioning that if the assumption {\bf(H1)}  holds, then $\sum_d\not=\emptyset$.

Recalling the estimate (\ref{2.2}), as the consequence of Lemmas \ref{l5.1}, \ref{l5.2} and \ref{l5.3}, we have

\begin{cor}\lbl{c5.1}\,{\rm(i)}\,If $h_\infty<\infty$, then $h_\infty\le h^*$. Hence, $h_0\ge h^*$ implies $h_\infty=\infty$ for all $\mu>0$;

{\rm(ii)}\,When $h_0<h^*$, there exist $\mu_0,\,\mu^0>0$, such that $h_\infty\leq h^*$ for $\mu\leq\mu_0$, and $h_\infty=\infty$ for $\mu\geq\mu^0$.
\end{cor}

Finally, we give the sharp criteria for spreading and vanishing.

\begin{theo}\lbl{th5.1}\,{\rm(i)}\, If $h_0\ge h^*$, then $h_\infty=\infty$ for all $\mu>0$;

{\rm(ii)}\, If $h_0<h^*$, then there exists $\mu^*>0$, such that $h_\infty=\infty$ for $\mu>\mu^*$, while $h_\infty\leq h^*$ for $\mu\leq\mu^*$.
\end{theo}

{\bf Proof}.\, Noticing Corollary \ref{c5.1}, by use of Lemma \ref{l3.1} and the continuity method, we can prove Theorem \ref{th5.1}. Please refer to the proof of Theorem 3.9 in \cite{DLin} for details. \ \ \ \fbox{}

\vskip 2pt Now we fix $h_0$, and regard $d$ and $\mu$ as the variable parameters to  describe the sharp criteria for spreading and vanishing. First of all, $\sum_{h_0}^+\not=\emptyset$ by Proposition \ref{p3.1}.

\begin{theo}\lbl{th5.2}\,{\rm(i)}\, When $d\in\sum_{h_0}^-$, we have $h_\infty=\infty$ for all $\mu>0$;

{\rm(ii)} For any fixed $d\in\sum_{h_0}^+$, there exists $\mu_0=\mu_0(d)>0$ such that $h_\infty<\infty$ provided $0<\mu\le\mu_0$. If, in addition, $\sum_d\not=\emptyset$ for such $d$, then there exists $\mu^*>0$, such that $h_\infty=\infty$ when $\mu>\mu^*$, and $h_\infty<\infty$ when $\mu\leq\mu^*$.
  \end{theo}

\begin{remark}\lbl{r5.1}\,{\rm(i)}\,Let $\hat x$ be given in {\bf(H2)}. From the proof of Proposition {\rm\ref{p3.3}} we see that if $h_0>\hat x$ then $\sum_{h_0}^-\not=\emptyset$, and there exists $d_0>0$ such that $d\in\sum_{h_0}^-$ for all $0<d\le d_0$.

{\rm(ii)}\, If the condition {\bf(H1)} holds, then $\sum_d\not=\emptyset$ for all $d>0$ by Proposition {\rm\ref{p3.2}}.
 \end{remark}

{\bf Proof of Theorem \ref{th5.2}}. (i)\, For any $d\in\sum_{h_0}^-$, we have $\lambda_1(h_0;d,a)\le 0$. If $\lambda_1(h_0;d,a)<0$, then $\sum_d\not=\emptyset$ and $h_0>h^*(d)$. If $\lambda_1(h_0;d,a)=0$, then $h_0=h^*(d)$. By Theorem \ref{th5.1}(i), $h_\infty=\infty$ for all $\mu>0$.

\vskip 2pt (ii)\, For the fixed $d\in\sum_{h_0}^+$, we have $\lambda_1(h_0;d,a)>0$. By Lemma \ref{l5.2}, there exists $\mu_0>0$ such that $h_\infty<\infty$ for $\mu\leq\mu_0$. If, in addition, $\sum_d\not=\emptyset$, then there exists $H\gg 1$ such that $\lambda_1(H;d,a)<0$. Take advantage of Lemma \ref{l5.3}, there exists $\mu^0>0$ such that $h_\infty>H$ provided $\mu\ge\mu^0$, which implies $\lambda_1(h_\infty;d,a)<\lambda_1(H;d,a)<0$. Hence, $h_\infty=\infty$ for $\mu\ge\mu^0$ by Lemma \ref{l5.1}. The remaining proof is the same as that of Theorem 3.9 in \cite{DLin}. \ \ \ \fbox{}

\vskip 2pt Now we fix $d$ and consider $h_0, \mu$ as the varying parameters. Note that (\ref{4.3}) implies {\bf(H1)}, combining Theorems \ref{th4.1}, \ref{th4.3} and \ref{th5.1}, we immediately obtain the following spreading-vanishing dichotomy and sharp criteria for spreading and vanishing.

\begin{theo}\lbl{th5.3} Assume that {\rm(\ref{4.3})} holds. Let $h^*$ be the unique positive root of $\lambda_1(\ell;d,a)=0$ and $(u,h)$ be the unique solution of {\rm(\ref{1.1})}. Then the following alternative holds:

Either\vspace{-2mm}\begin{itemize}
\item[{\rm(i)}]Spreading: $h_\infty=\infty$ and $\lim_{n\to\infty}u(t+nT,x)=U(t,x)$ uniformly on $[0,T]\times[0,L]$ for any $L>0$, where $U(t,x)$ is the unique $T$-periodic positive solution to {\rm(\ref{4.2})};\vspace{-2mm}\end{itemize}
or
\vspace{-2mm}\begin{itemize}
\item[{\rm(ii)}]Vanishing: $h_\infty\le h^*$ and $\lim_{t\to\infty}\,\max_{0\leq x\leq h(t)}u(t,x)=0$.
\end{itemize}\vspace{-2mm}

Moreover,

{\rm(iii)} If $h_0\ge h^*$, then $h_\infty=\infty$ for all $\mu>0$;

{\rm(iv)}\, If $h_0<h^*$, then there exist $\mu^*>0$, such that $h_\infty=\infty$ for $\mu>\mu^*$, while $h_\infty\leq h^*$ for $\mu\leq\mu^*$.
\end{theo}

When we fix $h_0$ and regard $d, \mu$ as the variable parameters. Note that (\ref{4.3}) implies {\bf(H1)}, using Theorems \ref{th4.1}, \ref{th4.3}, \ref{th5.2} and Remark \ref{r5.1}, we have the following spreading-vanishing dichotomy and sharp criteria.

\begin{theo}\lbl{th5.4}\,  Assume that $h_0>\hat x$, where $\hat x$ is given in {\bf(H2)}. Let the condition {\rm(\ref{4.3})} hold and $(u,h)$ be the unique solution of {\rm(\ref{1.1})}. Then the following alternative holds:

Either\vspace{-2mm}\begin{itemize}
\item[{\rm(i)}]Spreading: $h_\infty=\infty$ and $\lim_{n\to\infty}u(t+nT,x)=U(t,x)$ uniformly on $[0,T]\times[0,L]$ for any $L>0$, where $U(t,x)$ is the unique $T$-periodic positive solution to {\rm(\ref{4.2})};\vspace{-2mm}\end{itemize}
or
\vspace{-2mm}\begin{itemize}
\item[{\rm(ii)}]Vanishing: $h_\infty<\infty$ and $\lim_{t\to\infty}\,\max_{0\leq x\leq h(t)}u(t,x)=0$.
\end{itemize}\vspace{-2mm}

Moreover,

{\rm(iii)}\, If $d\in\sum_{h_0}^-$, then $h_\infty=\infty$ for all $\mu>0$;

{\rm(iv)} If $d\in\sum_{h_0}^+$, then there exists $\mu^*>0$, such that $h_\infty=\infty$ when $\mu>\mu^*$, and $h_\infty<\infty$ when $\mu\leq\mu^*$.
\end{theo}

In the last part of this section, we estimate the asymptotic spreading speed of the free boundary $h(t)$ when spreading occurs. To this aim, let us first state a known result, which plays an important role in the study of asymptotic spreading speed. For an integrable $T$-periodic function $\xi(t)$, we define
   \[\overline \xi=\frac 1T\int_0^T\xi(t){\rm d}t.\]
\begin{prop}\lbl{p5.1}{\rm(\cite[Section 2]{DGP})}\, Let $d>0$ and $0<\nu<1$ be the given constants. Assume that $p,q\in C^\nu([0,T])$ are positive $T$-periodic functions and $k\in C^\nu([0,T])$ is a nonnegative $T$-periodic function. Then the problem
\bess\left\{\begin{array}{ll}
 w_t-dw_{xx}+k(t)w_x=p(t)w-q(t)w^2, \ \ &0\le t\le T, \ 0<x<\infty,\\[1mm]
 w(t,0)=0, &0\le t\le T,\\[1mm]
w(0,x)=w(T,x),&0\le x<\infty
  \end{array}\right.\eess
has a positive $T$-periodic solution $w^k(t,x)\in C^{1,2}([0,T]\times[0,\infty))$ if and only if $\overline k<2\sqrt{d\overline p}$, and such a solution is unique when it exists. Furthermore, the following hold:

{\rm(i)} $w^k_x(t,x)>0$ and $w^k(t,x)\to v(t)$ uniformly on $[0,T]$ as $x\to\infty$, where $v(t)$ is the
unique positive periodic solution of the problem
 \[v'=p(t)v-q(t)v^2, \ \ 0\le t\le T; \ \ \ v(0)=v(T);\]

{\rm(ii)} For any given nonnegative $T$-periodic function $m\in C^\nu([0,T])$ satisfying $\overline m<2\sqrt{d\overline p}$, the assumption $m\le,\,\not\equiv k$ implies $w^m_x(t,0)>w^{k}_x(t,0)$, $w^m(t,x)>w^{k}(t,x)$ for $0\le t\le T$ and $x>0$;

{\rm(iii)} For each $\mu>0$, there exists a unique positive $T$-periodic function $k_0(t)=k_0(\mu,p,q)(t)\in C^\nu([0,T])$ such that $\mu w^{k_0}_x(t,0)=k_0(t)$ in $[0,T]$,  and $0<\overline k_0<2\sqrt{d\overline p}$.
  \end{prop}

\begin{theo}\lbl{th5.5} Assume that {\rm(\ref{4.3})} holds with $\rho=0$. When the spreading occurs, i.e., $h_\infty=\infty$, we have
 \bess
 \liminf_{t\to\infty}\frac{h(t)}t\ge \frac 1T\int_0^Tk_0(\mu,a_\infty,b^\infty)(t){\rm d}t, \ \ \ \limsup_{t\to\infty}\frac{h(t)}t\le
 \frac 1T\int_0^Tk_0(\mu,a^\infty,b_\infty)(t){\rm d}t,
  \eess
 \end{theo}

{\bf Proof}.\, This proof can be done by the same manner as section 4 of \cite{DGP}. Because the length is too long, we omit the details. The interested readers can refer to that reference.

\section{Concluding remarks}
\setcounter{equation}{0}{\setlength\arraycolsep{2pt}

We comment finally on some points raised by the theoretical investigation. Firstly, from the above discussion we have seen that $\lambda_1(\infty; d,a):=\lim_{\ell\to\infty}\lambda_1(\ell;d,a)<0$ is an essential condition. This number is only characterized by the dispersal rate $d$ and intrinsic growth rate $a(t,x)$, is independent of the self-limitation coefficient $b(t,x)$, moving parameter $\mu$ and initial value $u_0(x)$.

The main conclusions of this paper can be briefly summarized as follows:

(I)\, If one of the following holds:

(i)\, the diffusion rate $d\in\sum_{h_0}^-$ ($h_0$, $a$ are fixed);

(ii)\, the intrinsic growth rate $a(t,x)$ is suitable ``positive" in the sense of $\lambda_1(h_0;d,a)\le 0$ ($h_0$, $d$ are fixed);

(iii)\, the initial habitat $h_0$ is suitable ``larger"  in the sense of $\lambda_1(h_0;d,a)\le 0$ ($d$, $a$ are fixed),\\
then the new or invasive species will successfully spread and survive in the new environment (maintain a positive density distribution), regardless of initial population size $u_0(x)$ and value of the moving parameter $\mu$.

\vskip 2pt (II)\, When the above situations do not appear, we can find a critical $\mu^*$ such that the species will spread successfully when $\mu>\mu*$, and that the species fails to establish and  will be extinct in a long run when $\mu<\mu*$. A better way to modulate the moving parameter $\mu$ is to control the surrounding environment.

Now we analyze the relations and differences of the free boundary problem (\ref{1.1}) and the corresponding Cauchy problem and fixed domain problem, and make a comparison between our results and those obtained for the corresponding Cauchy problem and initial-boundary value problem with fixed boundary.

By the same arguments as that of Theorem 5.4 in \cite{DG1}, we can prove that the solution $u(t,x)$ of (\ref{1.1}) converges to $v(t,x)$ as $\mu\to\infty$, where $v(t,x)$ is the solution of the initial-boundary value problem in the half space
\bes
 \left\{\begin{array}{lll}
 v_t-d v_{xx}=a(t,x)v-b(t,x)v^2, \ \ &t>0,\ \ x>0,\\[1mm]
 B[v](t,0)=0,\ \ &t\ge 0,\\[1mm]
 v(0,x)=v_0(x),& x\ge 0,
 \end{array}\right.
 \lbl{6.1}\ees
where $v_0(x)=u_0(x)$ for $0\le x\le h_0$ and $v_0(x)=0$ for $x>h_0$. This shows that (\ref{6.1}) is the limiting problem of (\ref{1.1}) as $\mu\to\infty$. If $\alpha=0$, the problem (\ref{6.1}) can be extended into the Cauchy problem
 \bes\left\{\begin{array}{lll}
 v_t-d v_{xx}=a(t,x)v-b(t,x)v^2, \ \ &t>0,\ \ x\in\mathbb{R},\\[1mm]
  v(0,x)=v_0(x),& x\in\mathbb{R}.
 \end{array}\right.
 \lbl{6.2}\ees
Solutions of both (\ref{6.1}) and (\ref{6.2}) are positive for all $x$ once $t$ is positive. This seems to be a defect because the movement should be finite for any species.

When $a$ and $b$ are positive constants, the differential equation in (\ref{6.2}) is
 \bes
 u_t-d u_{xx} =u(a-bu), \ \  t > 0, \ x\in\mathbb{R},
 \lbl{6.3}\ees
A great deal of previous mathematical investigation on the spreading of population has been based on the traveling wave fronts and asymptotic spreading speed of (\ref{6.3}), please refer to, for example  Aronson \& Weinberger \cite{Ar, AW} for the model  (\ref{6.3}), Lewis et al. \cite{LLW} and Liang \& Zhao \cite{LZ} for the more general models. The known result for (\ref{6.3}) predicts successful spreading and establishment of the new species with any nontrivial initial population $u(0,x)$, regardless of its initial size and supporting area. However, this is not supported by empirical evidences, which suggest, in the contrary, that success of spreading depends on the initial size of the population; for example, the introduction of several bird species from Europe to North America in the 1900s was successful only after many or several initial attempts (cf. \cite{SK, LHP}, where more examples can be found).

On the other hand, when $\mu=0$, our free boundary problem (\ref{1.1}) reduces to the following initial-boundary value problem with fixed boundary:
 \bes\left\{\begin{array}{ll}
 w_t-dw_{xx}=a(t,x)w-b(t,x)w^2, \ \ &t>0, \ \, 0<x<h_0,\\[1mm]
 B[w](t,0)=w(t,h_0)=0,&t>0,\\[1mm]
  w(0,x)=u_0(x),&0\le x\le h_0
  \end{array}\right.
  \lbl{6.4}\ees
These discussions indicate that the free boundary problem (\ref{1.1}) is in the intermediate state of (\ref{6.1}) and (\ref{6.4}).

It is well known that (see Theorems 24.2 and 28.1 of \cite{Hess}) if $\lambda_1(h_0;d,a)\ge 0$, then the solution $w$ of (\ref{6.4}) must extinct in the long run, i.e., $w(t,x)$ tends to zero as $t\to\infty$ no matter how large the initial data $u_0(x)$ is and how small the self-limitation coefficient $b(t,x)$ is. This seems do not tally with the natural phenomenon.

For our free boundary problem (\ref{1.1}), the conclusions stated in Theorems \ref{th5.3} and \ref{th5.4} show that

(i) if $\lambda_1(h_0;d,a)\le 0$, then the spreading happens and  $u(t,x)$ tends to a positive $T$-periodic function for any moving parameter $\mu>0$ regardless of sizes of the initial data $u_0(x)$ and self-limitation coefficient $b(t,x)$;

(ii) if $\lambda_1(h_0;d,a)>0$, we can find a critical value $\mu^*(u_0,b)$ of the moving parameter, such that the spreading happens and $u(t,x)$ tends to a positive $T$-periodic function when $\mu>\mu^*(u_0,b)$, while the vanishing occurs and $u(t,x)$ will extinct in the long run for $\mu\le\mu^*(u_0,b)$.

From our proof we can see that the number $\mu^*(u_0,b)$ is decreasing in $u_0(x)$ and increasing in $b(t,x)$. The above conclusion (ii) tells us that, in the case $\lambda_1(h_0;d,a)>0$, whether the species establishes itself successfully or not depends on the sizes of the initial data and self-limitation coefficient. This seems to match the reality better and is supported by numerous empirical evidences introduced above.

Our conclusions indicate that the above mentioned shortcomings of (\ref{6.1}), (\ref{6.3}) and (\ref{6.4}) will not appear if we use the corresponding free boundary problem to describe the spreading or persistence instead of the Cauchy problem or initial-boundary value problem with fixed boundary. Furthermore, unlike the Cauchy problem model in which the spreading front is represented by an unspecified level set of the solution, the free boundary model gives an exact location of the spreading front $x=h(t)$ for any given time. We also remark that the important feature of (\ref{6.3}), namely the spreading front invades at a linear rate in time, is retained by the free boundary model.

We believe that the cause of the above mentioned defects for problems (\ref{6.1}), (\ref{6.3}) and (\ref{6.4}) is that the speed of the expanding front is infinity (too fast) when we use (\ref{6.1}) and (\ref{6.3}), and is zero (too slow) when we use  (\ref{6.4}). In the model (\ref{1.1}), the speed of the expanding front may be more reasonable.

\vskip 4pt \noindent {\bf Acknowledgment:} The author would like to thank the anonymous referees for their helpful comments and suggestions.

\bibliographystyle{elsarticle-num}

\end{document}